\documentclass[12pt, twoside,a4paper]{article}
\pdfoutput=1
\usepackage[cp1251]{inputenc}
\usepackage{ifthen}
\usepackage{amsmath}
\usepackage{amsfonts}
\usepackage{latexsym}
\usepackage{amsthm}
\usepackage{amssymb}
\newcommand{\CopyName}{ V.\ M.\ Zhuravlov}
\newcommand{\NAME}{ V.\ M.\ Zhuravlov}
\newcommand{\Year}{2024}
\newcommand{\rightheadtext}{Associative ideals in monoids}
     \newcounter{chapter}
     \newcounter{artpage}[chapter]
     \newcommand{\vs}{\vspace{.1in}}
     \newcommand{\vsk}{\vspace{.2in}}
\makeatletter
     \renewcommand{\@evenhead}{\footnotesize \ifthenelse{\value{artpage}=0}
     {\hfil}{\thepage\hfil \textsc {\leftmark} \hfil } }
     \renewcommand{\@oddhead}{\footnotesize\ifthenelse{\value{artpage}=0}
     {\hfil}{\hfil \textsc \rightmark \hfil \thepage} }
     \newcommand{\logo}{\baselineskip2pc \hbox to\hsize{\hfil\copyright\,\footnotesize
     \CopyName, \Year}}
     \renewcommand{\@oddfoot}{\ifthenelse{\value{artpage}=0}{\logo
     \refstepcounter{artpage}} {\hfil\refstepcounter{artpage}}}
     \renewcommand{\@evenfoot}{\ifthenelse{\value{artpage}=0}{\logo
     \refstepcounter{artpage}} {\hfil\refstepcounter{artpage}}}
     \renewcommand{\section}{\@startsection{section}{1}{0pt}{3.5ex plus
     1ex minus .2ex}{2.3ex plus 2.ex}{\large\hfil\textsc}}
     \vspace{3pt}\vs

\vs
\newcommand{\tit}{Associative ideals in monoids}
\date{2023}
\usepackage{hyperref}
\usepackage{graphicx}
\graphicspath{{Images/}}
\begin{document}
\hfill
\vspace{0.3in}
\markboth{{\NAME}}{{\rightheadtext}}\begin{center} \textsc {\CopyName} \end{center}\begin{center} \renewcommand{\baselinestretch}{1.3}\bf {\tit} \end{center}
\vspace{20pt plus 0.5pt} {\abstract{\noindent
The article investigates the properties of associative ideals in monoids. Such ideals have some applications in the logic of non-standard sequences and category theory. The relations of these ideals with the verbal structure of words over generators are examined. Finitely and infinitely generated monoids are treated separately.\newline
\textit{Associative ideals in monoids, 2024, msc: 20M12;\vspace{3pt}}\newline
}
}\vsk
\section{Introduction}\par
This article is a supplement to my article: \href{https://arxiv.org/abs/2312.00831}{"Predicates and terms from non-standard sequences"}. However, this supplement also has its own significance. Here we will discuss the properties of associative ideals of monoids, the notion of which was introduced in the aforementioned article. This will help to better understand the relation of monoids with category theory, logic and sequences. First, we will recall the basic definitions. Group theory is logically equivalent to semigroup theory. A semigroup with a unit is a monoid; the introduction of a unit (which is the ultimate universal constant for the whole set of elements with multiplication; these elements, in turn, uniquely define the only possible group extension — the introduction of inverse elements) is an inessential extension of the theory. The set of elements of a monoid is called an ideal if it contains all the products of its elements with any other elements of the monoid.\par
\textbf{Definition 1. A two-sided ideal \textit{Q} is called associative if:}
$$\forall a(a\neq E)\forall b(b\neq E)\forall c(c\neq E):(abc\in Q)\Longrightarrow[(ab\in Q)\vee(bc\in Q)]$$
where \textbf{\textit{E}} is the unit of the monoid \textbf{\textit{M}}.
\section{Properties of associative ideals}\par
In the aforementioned article, we proved that every prime ideal is associative, and that the converse is false. There are also associativity properties that are specific to prime ideals.\par
Let us begin our investigation by establishing some "ancillary" properties of associative ideals. Their intersection (as well as the intersection of prime ideals) is not necessarily associative; however, arbitrary unions are always associative (the same holds true for prime ideals). Such properties generate certain topologies on the monoid; monoid multiplication is continuous in these topologies.\par
It is also evident that the \textbf{property of being an associative ideal is preserved under epimorphisms of the monoid.}\par
Here are some additional propositions that hold true for associative ideals:
Let the ideal Q of the monoid M be associative. And let: $x,y,z\not\in Q$. Then:
$$(xyz\in(M-Q))\Longleftrightarrow[xy\in(M-Q)]\wedge[yz\in(M-Q)]$$
— this is correct, as if the composition of three elements belongs to an associative set, then the composition of any one pair of elements in their sequential arrangement also belongs to it (therefore, the converse implication is also true—derived from the negations of these statements).\par
\textbf{The essence of this proposition is that the complement to an associative ideal is an associative set, and this is "strong" associativity (in the sense that \textit{both} pairwise products are included in it); however, this is generally a strongly simple set:}
$$(xy\in (M-Q))\Longrightarrow(x\in (M-Q))\wedge(y\in (M-Q))$$
\textbf{— since \textit{Q} is a two-sided ideal.}\par
Logically, this signifies (refer again to the aforementioned article:\newline \href{https://arxiv.org/pdf/2312.00831.pdf}{arxiv.org/pdf/2312.00831.pdf}), that the existence of sequences in the left part of the equivalence implies the existence of sequences in the right part. A direct consequence of this proposition is:\par
\textbf{Assume now that the ideals} $Q_{1}\subseteq Q_{2}$ \textbf{are associative. Then the set} $(Q_{2}-Q_{1})$ \textbf{will also be associative.}\par
Let's proceed to the principal assertions of this article.\par
\textbf{The Boolean monoid of a monoid \textit{M} (being a Boolean algebra) is itself a monoid with the identity element \textit{E} and multiplication of two sets defined as the set of all products of elements from these sets. Moreover, the set of all ideals is a submonoid of this Boolean monoid.}\par
This is a rather self-evident statement.
\section{The connection of associative and prime ideals with the Cartesian product of the set of generators}\par
By the set of generators \textbf{\textit{S}} of the monoid \textbf{\textit{M}}, we mean any collection of elements from \textbf{\textit{M}} such that every element of M is a certain word formed from the alphabet \textbf{\textit{S}}. That is, we do not impose any conditions of minimality on \textbf{\textit{S}}. The set of generators for an ideal is the set of its elements such that the ideal is the union of their main ideals.\par
A generator will be called atomic if it has no divisors other than itself and the identity element \textbf{\textit{E}}.\par
Suppose the monoid M is generated by a finite or infinite set of atomic generators \textbf{\textit{S}}. Let's select a certain subset $Z\subseteq S$ of its generators. Then:\newline
\textbf{Lemma 1:} \textbf{\textit{The ideal closure I(Z) is a prime ideal. The ideal closure}} $I(Z\cdot Z)$ \textbf{\textit{is associative but not a prime ideal.}}\par
This has been demonstrated by case analysis of the inclusion of products from \textbf{\textit{Z}} in arbitrary words. Additionally, one can select a certain subset $F\subseteq(S \times S)$ and construct a "multiplication table":
$$s_1 \circ s_2 \equiv \chi_F (s_1, s_2)$$
— where $\chi_F$ is the characteristic function of \textbf{\textit{F}}. This leads to a certain generalization of Lemma 1:\par
\textbf{Theorem 1: \textit{For any monoid M that has a certain set of atomic generators S, and for any subset $F \subseteq (S \times S)$ — the collection of monoid products:}}
$$\{s_j s_k \mid ((s_j, s_k) \in F)\land(\chi_F (s_j, s_k) = 1)\}$$
\textbf{textit{constitutes the set of generators of a certain associative ideal Q}} (where the indexing is absolutely arbitrary—it does not have to be numerical or even countable.)\newline
Here, the elements of \textbf{\textit{M}} are represented by words composed of the generators \textbf{\textit{S}}. These words may be connected by certain relations — \textbf{\textit{M}} is not necessarily free. But if a word \textbf{\textit{x}} is expressible (even if not in a unique way) such that it includes the product $(s_j s_k)$ with the unit characteristic function \textbf{\textit{F}}, then $x \in Q$.\newline
Proof: The monoid \textbf{\textit{M}} is not necessarily free. However, each of its elements equals some word from \textbf{\textit{S}}. Suppose $xyz \in Q$, where $Q = I(F)$ — the ideal closure of a subset of pairs of generators. Then there exists a representation of \textbf{\textit{xyz}} as a word containing the product of a certain pair $s_j s_k$ from \textbf{\textit{F}}. Each of the elements \textbf{\textit{(x, y, z)}} will also be a word written in the alphabet \textbf{\textit{S}}. Where in these words is $s_j s_k$ located? If in \textbf{\textit{x}}, then $xy \in Q$, since \textbf{textit{Q}} is a two-sided ideal. If in \textbf{\textit{z}}, then similarly $yz \in Q$. And if in \textbf{\textit{y}}, then $(xy \in Q) \land (yz \in Q)$. If $s_j$ and $s_k$ are included in different elements of the triple $(x, y, z)$ (while remaining adjacent in the word \textbf{\textit{(xyz)}}), then $(xy \in Q) \lor (yz \in Q)$. In any case, the ideal \textbf{\textit{Q}} will be associative. The proof is concluded.
\section{On infinite divisibility}\par
It is noteworthy that we are discussing a monoid with atomic generators. However, this also pertains to infinitely divisible collections of generators... The divisibility relation $(x\preceq y)\equiv\exists a:((xa = y)\vee(ax = y))$ is a pre-order relation. In such a defined pre-order, being an ideal means containing the right segment of each of its elements. Regarding the left segment of each element, the simplicity of the ideal implies the divisibility of each element within this ideal, i.e., the maximality of the ideal relative to each of its elements. The associativity of an ideal (as a subsequent, weaker gradation of divisibility) signifies the next, weaker (after maximal) gradation in the size of the ideal. If there is an infinite set of atomic generators of the monoid \textbf{\textit{M}} (incomparable in the aforementioned pre-order), then the theorem remains true for any (finite or infinite) \textbf{\textit{S}}. However, if \textbf{\textit{S}} contains infinitely descending chains, then these chains will (or will not) have atomic limits (which, in certain cases, depends on whether we accept the truth or falsity of Zorn's Lemma for \textbf{\textit{S}}). If such limits exist, then they can be used to restrict the entire set \textbf{\textit{S}}, and our theorem remains true. Otherwise, it is necessary to introduce for \textbf{\textit{F}} conditions similar to the Cauchy conditions that define limits in mathematical analysis. Namely, \textbf{\textit{F}} must be such a (naturally, infinite) subset that:
$$\forall s_{j}\forall (s_{k}:(s_{j},s_{k})\in F)\Longrightarrow\exists t_{j}\exists t_{k}:((t_{j},t_{k})\in F)\wedge(t_{j}\prec s_{j})\wedge(t_{k}\prec s_{k}))$$
— that is, along with every pair of generators \textbf{\textit{F}}, it must also contain a pair of strictly smaller (and consequently, not equal to them) generators, in the pre-order specified above. When this condition is satisfied, the theorem remains true even in the infinite case.\par
And then we encounter a situation completely analogous to Theorem 1. Therefore, the following theorem is valid:\par
\textbf{Theorem 2:}\textbf{\textit{ An ideal is prime if and only if its generators are a subset of the monoid's generators and one of the following three conditions is met:}}\par
\textbf{\textit{1. Along with each of its generators, it contains a strictly decreasing (finite or infinite) chain of the monoid's generators, the supremum of which is the original generator of the ideal;}}\par
\textbf{\textit{2. Or the generator of the ideal is an atom;}}\par
\textbf{\textit{3. Or the chain has an atomic limit that belongs to the ideal.}}\par
\textbf{\textit{An ideal is associative if and only if its generators are pairwise products of some subset of such generators.}}\par
Proof: Let \textbf{\textit{P}} be an ideal, and let \textbf{\textit{S}} be the (finite or infinite) set of generators of the monoid \textbf{\textit{M}}. Then the primicity of \textbf{\textit{P}} is equivalent to the following condition: The generators \textbf{\textit{s}} of the ideal \textbf{\textit{P}} are a subset of \textbf{\textit{S}}, and along with each such generator, \textbf{\textit{P}} contains an infinitely descending chain from \textbf{\textit{S}}, the supremum of which is \textbf{\textit{s}}; if no such chains exist, then \textbf{\textit{s}} is an indivisible atom. This is by definition of a prime ideal.\par
Associative ideals have generators that are pairwise products of some of these chains (or possibly atomic elements of \textbf{\textit{S}}).\par
Let's articulate this more precisely. Denote by $C_{s}$ an infinite, strictly descending chain of generators in the preorder of divisibility, the supremum of which is $s \in C_{s}$ (if such a chain exists, it may not be unique). If such a chain has a lower bound $s_{o}$, which is an atomic generator, we replace the element \textbf{\textit{s}} in \textbf{\textit{S}} with $s_{o}$ — as the generator of the ideal.\par
If no such lower bound exists, there are two paths to follow. We can introduce it artificially as a non-standard element and replace the entire chain with it (along with \textbf{\textit{s}}). Alternatively, we can expand the set of generators of the ideal by including the entire chain. As a result, we obtain the set \textbf{\textit{S0}}, which will also be the subset of generators for \textbf{\textit{M}} and will be a set of generators for the ideal. The proof is concluded.\par
As we can see, simple and associative ideals are closely related — the sets of their generators are sufficiently "close" to each other.
\section{Roots of monoids and quadratic monoids}\par
In this context, ideas for some new directions in the study of monoid ideals emerge.\newline
\textbf{Definition 2:\textit{ The root}} $\sqrt{M}$\textit{\textbf{ of a monoid M is defined as the smallest monoid that is free relative to M — such that each generator of M is a pairwise product of the generators of}} $\sqrt{M}$\textit{\textbf{. Similarly, the smallest monoid}} $M^2$\textit{\textbf{, which is free relative to M and whose generators are all pairwise products of the generators of M, is called a quadratic monoid.}} The term "free relative to \textit{\textbf{M}}" implies that $\sqrt{M}$ has no relations other than those present in \textit{\textbf{M}} (and, naturally, "pairwise" relations between the generators of the monoids).\par
It is clear that $M^2 \subseteq M \subseteq \sqrt{M}$. Now consider a prime ideal \textbf{\textit{P}} in \textbf{\textit{M}}. Corresponding to it are the ideals: $\sqrt{P}$ — the minimal ideal in $\sqrt{M}$ containing \textbf{\textit{P}}; and $P^2$ — the principal ideal in $M^2$ for the set of pairwise products of the generators of \textbf{\textit{P}}. The following corollary from the theorem appears evident:\newline
\textbf{Corollary:} $\sqrt{P}$ \textbf{\textit{will be associative if and only if P is prime. Conversely, P will be associative if and only if}} $P^2$ \textbf{\textit{is prime.}}
\section{Logical Notes}\par
We have already devised a scheme for substituting constants with a pair of functions (right and left) to localize the zero term in non-standard sequences (see article cited above). Moreover, we understand that all such sequences boil down to binary relations between elements. Additionally, we can model sequences not only with a monoid but also with a semigroup (an associative operation without an identity element — thus, the identity mapping is not considered when modeling mappings with category arrows). Nonetheless, we cannot escape certain limiting relations. Specifically, let the following axiom hold true in the standard model:
$$\forall x\forall y:A(x,y)$$
It is replaced by:
$$\forall x\forall y:e(x,y)\Longrightarrow A(x,y)$$
— where e(x,y) is an existence predicate; but this is not enough. We must also impose a condition on this predicate:
$$e(x)\Longrightarrow\exists y\exists z:e(y,x)\wedge e(x,z)$$
That is, any existing sequence can be extended to the left and to the right. Consequently, for any \textbf{\textit{x}} and \textbf{\textit{y}}, there exist such paired elements for which \textbf{\textit{A(x, y)}} is indeed satisfied. This also applies to formulas from longer sequences.\par
There are also unary terms that emerge in place of 0-ary ones; for these, too, there always exist sequences with those elements, the functions of which they represent.\par
In general, let there be an axiom:
$$\bot x_{j|j\in(1,...,n)}:A(x_{j|j\in(1,...,n)})$$
where: $\bot x_{j}$ means an arbitrary quantifier of the variable $x_{j}$;\par
We also see that the pairwise products of any subset of a given set of monoid generators can serve as generators for a certain associative ideal; thus, complementing this subset to the entire set will create an associative ideal, in some sense complementary to another. Moreover, a "complete" multiplication table of the entire given set of generators can also create associative ideals.\par
These are, in my view, the fundamental properties of associative ideals of monoids; and possibly, these properties will find application in category theory and local sequences in general.
\section{References}
\noindent [1]   Hanna Neumann, Springer-Verlag,\newline
\textit{ Varieties of Groups}, (1967).\newline
[2] Serge Lang, Springer-Verlag,\newline
\textit{Algebra}, (1965).\newline
[3] P.M. Kon,  D. Reidel Publishing Company,\newline
 \textit{Universal Algebra}, (1961).\newline
\end{document}